\documentclass[12pt, a4paper, reqno]{amsart}
\usepackage{amsfonts,amsmath,latexsym,verbatim,amscd,mathrsfs,color,array}
\usepackage[colorlinks=true]{hyperref}
\usepackage{amssymb,amsfonts,graphicx}
\usepackage[hmargin=2.5cm, vmargin=2.5cm]{geometry}
\usepackage{bbm}
\usepackage{pdfsync}
\usepackage{epstopdf}
\usepackage{cite}
\usepackage{multirow}
\usepackage[font=bf,aboveskip=15pt]{caption}
\usepackage[toc,page]{appendix}
\allowdisplaybreaks[3]

\usepackage{hyperref}
\hypersetup{
    colorlinks=true,
    linkcolor=red,
    citecolor=blue,
    urlcolor=black
}

\usepackage{titlesec}
\titleformat{\section}
  {\large\bfseries\raggedright} 
  {\thesection.} 
  {0.4em} 
  {} 

\usepackage{fancyhdr}

\newtheorem{theo}{Theorem}[section]

\newtheorem{rema}{Remark}[section]

\newtheorem{defi}{Definition}[section]

\numberwithin{equation}{section}

\usepackage{marvosym} 

\usepackage{xcolor} 

\begin{document}

\title{\large Complete Classification of the Symmetry Groups of Monge-Amp\`{e}re Equation and Affine Maximal type Equation}
\author{Huan-Jie Chen$^1$ and Shi-Zhong Du$^1$}

\thanks{\Letter\enspace Shi-Zhong Du (szdu@stu.edu.cn) \enspace\&\enspace Huan-Jie Chen (19hjchen@stu.edu.cn)}
\thanks{$^1$\enspace\enspace The Department of Mathematics, Shantou University, Shantou 515063, P. R. China.}
\thanks{*\enspace\enspace This work is partially supported by Natural Science Foundation of China (12171299)}

\maketitle

\renewenvironment{abstract}
 {\small\noindent\textbf{Abstract}\enspace}
 {\par\vspace{0.9em}}

\begin{abstract}
  The affine maximal type hypersurface has been a core topic in Affine Geometry. When the hypersurface is presented as a regular graph of a convex function $u$, the statement that the graph is of affine maximal type is equivalent to the statement that $u$ satisfies the fully nonlinear partial differential equation
  \begin{equation}\label{e0.1}
       D_{ij}(U^{ij}w)=0, \ \ w\equiv[\det D^2u]^{-\theta}, \ \ \theta>0, \ \ \forall x\in{\mathbb{R}}^N
  \end{equation}
  of fourth order. This equation can be regarded as a generalization of the $N$-dimensional Monge-Amp\`{e}re equation
  \begin{equation}\label{e0.2}
      \det D^2u=1, \ \ \forall x\in{\mathbb{R}}^N
  \end{equation}
  of second order, since each solution of \eqref{e0.2} satisfies \eqref{e0.1} automatically. In this paper, we will determine the symmetry groups of these two important fully nonlinear equations without asymptotic growth assumption. Our method develops the Lie's theory to fully nonlinear PDEs.
\end{abstract}


\vspace{9pt}


\vspace{12pt}

\section{Introduction}

\noindent The fourth order fully nonlinear equation of affine maximal hypersurfaces is given by
   \begin{equation}\label{e1.1}
       D_{ij}(U^{ij}w)=0, \ \ \forall x\in\Omega\subset{\mathbb{R}}^N,
   \end{equation}
where $U^{ij}$ is the co-factor matrix of $u_{ij}$ and $w\equiv[\det D^2u]^{-\theta}, \theta>0$. Equation \eqref{e1.1} is the Euler-Lagrange equation of the affine area functional
   \begin{eqnarray*}
     {\mathcal{S}}(u,\Omega)&\equiv&\int_\Omega[\det D^2u]^{1-\theta}\\
      &=&\int_{{\mathcal{M}}_\Omega}K_0^{1-\theta}(1+|Du|^2)^{\vartheta}dV_{g_0},\ \ \vartheta=\frac{N+1}{2}-\frac{N+2}{2}\theta
   \end{eqnarray*}
for $\theta\not=1$ and
    $$
     {\mathcal{S}}(u,\Omega)\equiv\int_\Omega\log\det D^2u
    $$
for $\theta=1$, where $g_0$ and $K_0$ are the induced Riemannian metric and the Gauss curvature of the graph ${\mathcal{M}}_\Omega\equiv\big\{(x,z)\in{\mathbb{R}}^{N+1}: z=u(x), x\in\Omega\big\}$ respectively. Noting that
   $$
     D_jU^{ij}=0, \ \ \forall i=1,2,\cdots, N,
   $$
we can rewrite equation \eqref{e1.1} as
$$U^{ij}D_{ij}w=0,$$
or equivalently
   \begin{equation}\label{e1.2}
      u^{ij}D_{ij}w=0
   \end{equation}
for $[u^{ij}]$ representing the inverse of Calabi's metric $[u_{ij}]$ of graph ${\mathcal{M}}_\Omega$.

The special affine maximal case $\theta\equiv\frac{N+1}{N+2}$ has been extensively studied in the past. If one introduces the affine metric
  $$
    A_{ij}=\frac{u_{ij}}{[\det D^2u]^{1/(N+2)}}
  $$
on ${\mathcal{M}}_\Omega$ comparing to the Calabi's metric $g_{ij}=u_{ij}$ and sets
  $$
    H\equiv[\det D^2u]^{-1/(N+2)},
  $$
then it follows that \eqref{e1.2} reduces to
   \begin{equation}
    \triangle_{{\mathcal{M}}}H=0
   \end{equation}
for Laplace-Beltrami operator
   \begin{eqnarray*}
     \triangle_{{\mathcal{M}}}&\equiv&\frac{1}{\sqrt{A}}D_i(\sqrt{A}A^{ij}D_j) \\
     &=&HD_i(H^{-2}u^{ij}D_j)
   \end{eqnarray*}
with respect to the affine metric, where $A=\det A_{ij}$ and $[A^{ij}]$ stands for the inverse matrix of $[A_{ij}]$. So, the hypersurface ${\mathcal{M}}$ is affine maximal if and only if $H$ is harmonic on ${\mathcal{M}}_\Omega$.

The central topic in the study of the affine maximal type hypersurfaces is the Bernstein problem. A conjecture originally proposed by Chern \cite{Ch} for entire graph and then reformulated by Calabi \cite{Ca} to its fully generality asserts that any Euclidean complete, affine maximal type, locally uniformly convex $C^4$-hypersurface in ${\mathbb{R}}^{N+1}$ must be an elliptic paraboloid. This conjecture was completely resolved in the corner stone paper Trudinger-Wang \cite{TW1} for dimension $N=2$ with $\theta=3/4$, and later extended by Jia-Li \cite{JL2} to $N=2$ and $\theta\in(3/4,1]$ (see also Zhou \cite{Z} for a different proof). Over the past two decades, much efforts were done toward higher dimensional issues but not really successful yet, even for the case of dimension $N=3$. Partial affirmative results in this conjecture can be found in \cite{A,Do,JL,TW2,M}. Conversely, based on the priori fact that $SO(N)$ is one of the symmetry groups of geometric equation \eqref{e1.1}, the second author \cite{Du1,Du2} gave an opposite answer to Bernstein problem for the range $N\geq 2,\theta\in(0,(N-1)/N]$ by constructing various Euclidean complete group invariant solutions that are not elliptic paraboloid.

The affine maximal type equation can be regarded as a generalization of the Monge-Amp\`{e}re equation
   \begin{equation}\label{e1.4}
      \det D^2u=1, \ \ \forall x\in\Omega\subset{\mathbb{R}}^N
   \end{equation}
of second order, since each solution to \eqref{e1.4} satisfies the equation \eqref{e1.2}. There are many works about this classical second order partial differential equation, see for examples \cite{CNS,CW1,L3,TW3,P,T,Tso,W} about the existence and regularity of Dirichlet problem or Neumann problem, as well as the relevant therein.

On the other orientation, a classical method to find symmetry reductions of partial differential equations is the Lie group method. This method has been extensively investigated and widely applied to differential equations on Euclidean space as well as to geometric equations on Riemannian manifold. For example, Wo-Yang-Wang \cite{WYW} obtained group invariant solutions for centro-affine invariant flow by classifying symmetry groups and optimal systems within the framework of Lie's theory. The readers may also refer to other papers \cite{CL,CQ,CM,CDG,OST} and the elegant book \cite{O} for the related topic. In the frame of symmetry group, Lie's theory illustrates all transformations preserving the equations invariant, and then classifies the corresponding group invariant solutions. Although there are vast literatures about the affine maximal type equation and the Monge-Amp\`{e}re equation, there are few papers discuss the symmetry groups of them using the Lie's theory. This is the main purpose of this paper to classify the symmetry groups of \eqref{e1.2} and \eqref{e1.4} completely.

For the Monge-Amp\`{e}re equation, the following classification result is shown.
\begin{theo}\label{t1.1}
  The symmetry group of the Monge-Amp\`{e}re equation
    \begin{equation}\nonumber
      \det D^2u=1, \ \ \forall x\in\mathbb{R}^N,
    \end{equation}
  are generated by
   \begin{eqnarray*}
     g^i_1(\varepsilon) &: & (x,u)\to(\widetilde{x},\widetilde{u})=(x^1, \cdots, x^i+\varepsilon, \cdots, x^N, u), \\
     g_2(\varepsilon) &: & (x,u)\to(\widetilde{x},\widetilde{u})=(x, u+\varepsilon),\\
     g^i_3(\varepsilon) &: & (x,u)\to(\widetilde{x},\widetilde{u})=(x, u+\varepsilon x^i),\\
     g_4(\varepsilon) &: & (x,u)\to(\widetilde{x},\widetilde{u})=(A_\varepsilon x, u), \ A_\varepsilon\in SL(N),\\
     g^i_5(\varepsilon) &: & (x,u)\to(\widetilde{x},\widetilde{u})=(x^1, \cdots, e^{N\varepsilon}x^i, \cdots, x^N, e^{2\varepsilon}u),
   \end{eqnarray*}
  where $i=1,2,\cdots,N$.
\end{theo}

For the affine maximal type equation, the following classification result is shown.
\begin{theo}\label{t1.2}
  The symmetry group of the affine maximal type equation
    \begin{equation}\nonumber
      u^{ij}D_{ij}w=0, \ \ \forall x\in{\mathbb{R}}^N,
    \end{equation}
  are generated by
   \begin{eqnarray*}
     g^i_1(\varepsilon) &: & (x,u)\to(\widetilde{x},\widetilde{u})=(x^1, \cdots, x^i+\varepsilon, \cdots, x^N, u),\\
     g_2(\varepsilon) &: & (x,u)\to(\widetilde{x},\widetilde{u})=(x, u+\varepsilon),\\
     g_3(\varepsilon) &: & (x,u)\to(\widetilde{x},\widetilde{u})=(x, e^\varepsilon u),\\
     g^i_4(\varepsilon) &: & (x,u)\to(\widetilde{x},\widetilde{u})=(x, u+\varepsilon x^i), \\
     g_5(\varepsilon) &: & (x,u)\to(\widetilde{x},\widetilde{u})=(B_\varepsilon x, u),\ B_\varepsilon\in GL(N),\\
     g_6(\varepsilon) &: & (x,u)\to(\widetilde{x},\widetilde{u})=C_\varepsilon(x, u), \ C_\varepsilon\in SO(N+1)
   \end{eqnarray*}
 for $\theta=\frac{N+1}{N+2}$, and $g^i_1(\varepsilon)$-$g_5(\varepsilon)$ for $\theta\neq\frac{N+1}{N+2}$, where $i=1,2,\cdots,N$.
\end{theo}

Since no asymptotic growth assumption is imposed, our classification results are complete. The contents of the paper are organized as follows. First, we will recall some fundamental facts about Lie's theory for partial differential equations in Section \ref{sec2}. Second, we will calculate the coefficients from the fourth order prolongation formula in Section \ref{sec3}. Third, we will classify the symmetry groups of Monge-Amp\`{e}re equation and affine maximal type equation, thereby completing the proofs of Theorems \ref{t1.1}-\ref{t1.2} in Sections \ref{sec4}-\ref{sec5} respectively.

\vspace{20pt}

\section{Preliminary facts to Lie's theory on PDEs}\label{sec2}

\noindent  Let us recall some fundamental facts of Lie's theory to partial differential equations. We suggest the reader refer to the classical book \cite{O} of Peter Olver for the detailed introduction on the symmetry group of PDEs. At the beginning, we present the definition of a symmetry group for the equation as follows.
\begin{defi}
Let $F: {\mathbb{R}}^n\times{\mathbb{R}}^m\times{\mathbb{R}}^{mn}\times{\mathbb{R}}^{mn^2}\times\cdots\times{\mathbb{R}}^{mn^k} \to{\mathbb{R}}^l$ be a mapping from ${\mathbb{R}}^{n+m\frac{1-n^{k+1}}{1-n}}$ to ${\mathbb{R}}^l$. A symmetry group of the equation
$$F(x, u, Du, \cdots, D^ku)=0$$
is a local group of transformations $G$ defined on an open subset $M\subset\big\{(x,u)\in\mathbb{R}^n\times\mathbb{R}^m\big\}$ of the space of independent and dependent variables for the equation with the property that if $u=f(x)$ is a solution to the equation $F=0$, then whenever $g\circ f$ is defined for $g\in G$, we have $u=g\circ f(x)$ is also a solution.
\end{defi}

If one considers a partial differential equation
   \begin{equation}\label{e2.1}
      F(x, u, Du, \cdots, D^ku)=0
   \end{equation}
of $k$-order and supposes that
   \begin{equation}\label{e2.2}
    \overrightarrow{v}=\xi^i(x,u)\frac{\partial}{\partial x^i}+\phi(x,u)\frac{\partial}{\partial u}
   \end{equation}
is the infinitesimal generator of one-parameter group action $g(\varepsilon)\ (\varepsilon\in\mathbb{R})$ of \eqref{e2.1}, we have the $k$-th prolongation of $\overrightarrow{v}$ is the vector field
   \begin{equation}\label{e2.3}
       pr^{(k)}\overrightarrow{v}=\overrightarrow{v}+\phi^J(x,u^{(k)})\frac{\partial}{\partial u_J}
   \end{equation}
by prolongation formula in \cite{O} (Theorem $2.36$, Page $110$), where
    \begin{equation}\label{e2.4}
       \phi^J(x,u^{(k)})=D_J(\phi-\xi^iu_i)+\xi^iD_iu_J
    \end{equation}
for multi-indices $J=(j_1,\cdots, j_s), \ 1\leq j_s\leq n, 1\leq s\leq k$ as usually. Moreover, $g(\cdot)$ is an one-parameter symmetry group of the equation \eqref{e2.1} if and only if
   \begin{equation}
      pr^{(k)}\overrightarrow{v}F(x,D^ku)=0
   \end{equation}
holds for any $u^{(k)}\equiv(u, Du, \cdots, D^ku)$ satisfying
   \begin{equation*}
     F(x, u^{(k)})=0,
   \end{equation*}
where $x, u, Du, \cdots, D^ku$ are regarded as independent variables.

\vspace{20pt}

\section{Fourth order prolongation formula}\label{sec3}

\noindent To classify the symmetry groups of Monge-Amp\`{e}re equation and affine maximal type equation, we need to impose the second order and fourth order prolongation formulas respectively. In this section, we present only the fourth order prolongation formula, as it coincides with the second order one.

The $4$-th prolongation of the infinitesimal generator \eqref{e2.2} is expressed as
   $$
    pr^{(4)}\overrightarrow{v}=\xi^i\frac{\partial}{\partial x^i}+\phi\frac{\partial}{\partial u}+\phi^i\frac{\partial}{\partial u_i}+\phi^{ij}\frac{\partial}{\partial u_{ij}}+\phi^{ijk}\frac{\partial}{\partial u_{ijk}}+\phi^{ijkl}\frac{\partial}{\partial u_{ijkl}}.
   $$
By \eqref{e2.3}, the coefficients $\phi^i$ and $\phi^{ij}$ are determined by
    \begin{eqnarray}\nonumber
       \phi^i&=&D_i(\phi-\xi^su_s)+\xi^su_{is}\\
       &=&\phi_i+\phi_uu_i-(\xi^s_i+\xi^s_uu_i)u_s,\\ \nonumber
       \phi^{ij}&=&D_{ij}(\phi-\xi^su_s)+\xi^su_{ijs}\\ \nonumber
       &=&\phi_{ij}+\phi_{uj}u_i+\phi_{iu}u_j+\phi_{uu}u_iu_j+\phi_uu_{ij}\\
       &&-(\xi^s_{ij}+\xi^s_{uj}u_i+\xi^s_{iu}u_j+\xi^s_{uu}u_iu_j+\xi^s_uu_{ij})u_s\\ \nonumber
       &&-(\xi^s_i+\xi^s_uu_i)u_{js}-(\xi^s_j+\xi^s_uu_j)u_{is}.
    \end{eqnarray}
Since the expressions for the coefficients $\phi^{ijk}$ and $\phi^{ijkl}$ are relatively complicated, we need to introduce some notations for brevity. Let $\textstyle\sum_{\mathrm{circ}}\mathcal{A}(i,j,k)$ and $\textstyle\sum_{\mathrm{circ}}\mathcal{B}(i,j,k,l)$ represent the distinct circle summations over all permutations of indices $(i,j,k)$ and $(i,j,k,l)$ respectively. It is notable that in the permutations of $\phi_{ij}u_k$, there are totally six elements
   $$
    \phi_{ij}u_k, \ \ \phi_{jk}u_i, \ \ \phi_{ki}u_j,\ \ \phi_{ji}u_k, \ \ \phi_{ik}u_j, \ \ \phi_{kj}u_i.
   $$
Noting that $\phi_{ij}u_k=\phi_{ji}u_k, \ \ \phi_{jk}u_i=\phi_{kj}u_i, \ \ \phi_{ki}u_j=\phi_{ik}u_j$, the terms appeared in the circle summation $\sum_{circ}\phi_{ij}u_k$ will be only three terms. More precisely, we have the circle summations like
\begin{equation*}
\begin{aligned}
&\textstyle\sum\limits_{\mathrm{circ}}\phi_{ij}u_k=\phi_{ij}u_k+\phi_{ik}u_j+\phi_{jk}u_i,\\
&\textstyle\sum\limits_{\mathrm{circ}}\phi_{ij}u_{kl}=\phi_{ij}u_{kl}+\phi_{ik}u_{jl}+\phi_{il}u_{jk}+\phi_{jk}u_{il}+\phi_{jl}u_{ik}+\phi_{kl}u_{ij},\\
&\textstyle\sum\limits_{\mathrm{circ}}(u_{ij}u_{kl}+u_{ijkl})=(u_{ij}u_{kl}+u_{ik}u_{jl}+u_{il}u_{jk})+u_{ijkl}.
\end{aligned}
\end{equation*}

Now, the third coefficient $\phi^{ijk}$ can be written as
    \begin{eqnarray}\nonumber
      \phi^{ijk}&=&D_{ijk}(\phi-\xi^su_s)+\xi^su_{ijks}\\
       &=&\phi_{ijk}+\mathcal{T}_1+\mathcal{T}_2+\mathcal{T}_3,
    \end{eqnarray}
where the terms involving partial derivatives of $u$ are decomposed into three parts based on the highest derivative order:
\begin{equation*}
\left\{
\begin{aligned}
\mathcal{T}_1=&\textstyle\sum\limits_{\alpha=1}^2\textstyle\sum\limits_{\mathrm{circ}}\mathcal{A}^{(1)}_\alpha(i,j,k),\\
\mathcal{T}_2=&\textstyle\sum\limits_{\alpha=1}^2\textstyle\sum\limits_{\mathrm{circ}}\mathcal{A}^{(2)}_\alpha(i,j,k),\\
\mathcal{T}_3=&\textstyle\sum\limits_{\alpha=1}^2\textstyle\sum\limits_{\mathrm{circ}}\mathcal{A}^{(3)}_\alpha(i,j,k).
\end{aligned}
\right.
\end{equation*}
For the first order derivative term $\mathcal{T}_1$, its components are constituted by
\begin{eqnarray*}
\mathcal{A}^{(1)}_1(i,j,k)&=&(\phi_{jku}+\phi_{kuu}u_j)u_i-(\xi^s_{jku}+\xi^s_{kuu}u_j)u_su_i,\\
\mathcal{A}^{(1)}_2(i,j,k)&=&(\phi_{uuu}-\xi^s_{uuu}u_s)u_iu_ju_k-\xi^s_{ijk}u_s.
\end{eqnarray*}
For the second order derivative term $\mathcal{T}_2$, we have
\begin{eqnarray*}
\mathcal{A}^{(2)}_1(i,j,k)&=&(\phi_{ku}+\phi_{uu}u_k)u_{ij}-(\xi^s_{ku}+\xi^s_{uu}u_k)u_su_{ij},\\
\mathcal{A}^{(2)}_2(i,j,k)&=&-(\xi^s_{jk}+\xi^s_{ku}u_j+\xi^s_{uu}u_ju_k+\xi^s_uu_{jk})u_{is}.
\end{eqnarray*}
For the third order derivative term $\mathcal{T}_3$, we obtain that
\begin{eqnarray*}
\mathcal{A}^{(3)}_1(i,j,k)&=&(\phi_u-\xi^s_uu_s)u_{ijk},\\
\mathcal{A}^{(3)}_2(i,j,k)&=&-(\xi^s_i+\xi^s_uu_i)u_{jks}.
\end{eqnarray*}

Similarly, the fourth coefficient $\phi^{ijkl}$ is given by
    \begin{eqnarray}\nonumber
      \phi^{ijkl}&=& D_{ijkl}(\phi-\xi^su_s)+\xi^su_{ijkls}\\
      &=&\phi_{ijkl}+\mathcal{F}_1+\mathcal{F}_2+\mathcal{F}_3+\mathcal{F}_4,
    \end{eqnarray}
where the terms containing partial derivatives of $u$ are decomposed into four parts according to the highest order of the derivative:
\begin{equation*}
\left\{
\begin{aligned}
\mathcal{F}_1=&\textstyle\sum\limits_{\alpha=1}^3\textstyle\sum\limits_{\mathrm{circ}}\mathcal{B}^{(1)}_\alpha(i,j,k,l),\\
\mathcal{F}_2=&\textstyle\sum\limits_{\alpha=1}^4\textstyle\sum\limits_{\mathrm{circ}}\mathcal{B}^{(2)}_\alpha(i,j,k,l),\\
\mathcal{F}_3=&\textstyle\sum\limits_{\alpha=1}^2\textstyle\sum\limits_{\mathrm{circ}}\mathcal{B}^{(3)}_\alpha(i,j,k,l),\\
\mathcal{F}_4=&\textstyle\sum\limits_{\alpha=1}^2\textstyle\sum\limits_{\mathrm{circ}}\mathcal{B}^{(4)}_\alpha(i,j,k,l).
\end{aligned}
\right.
\end{equation*}
For the first order derivative term $\mathcal{F}_1$, the components are determined by
\begin{eqnarray*}
\mathcal{B}^{(1)}_1(i,j,k,l)&=&(\phi_{jklu}+\phi_{kluu}u_j+\phi_{luuu}u_ju_k)u_i,\\
\mathcal{B}^{(1)}_2(i,j,k,l)&=&-(\xi^s_{jklu}+\xi^s_{kluu}u_j+\xi^s_{luuu}u_ju_k)u_su_i,\\
\mathcal{B}^{(1)}_3(i,j,k,l)&=&(\phi_{uuuu}-\xi^s_{uuuu}u_s)u_iu_ju_ku_l-\xi^s_{ijkl}u_s.
\end{eqnarray*}
For the second order derivative term $\mathcal{F}_2$, we have
\begin{eqnarray*}
\mathcal{B}^{(2)}_1(i,j,k,l)&=&(\phi_{klu}+\phi_{luu}u_k+\phi_{uuu}u_ku_l+\phi_{uu}u_{kl})u_{ij},\\
\mathcal{B}^{(2)}_2(i,j,k,l)&=&-(\xi^s_{klu}+\xi^s_{luu}u_k+\xi^s_{uuu}u_ku_l+\xi^s_{uu}u_{kl})u_su_{ij},\\
\mathcal{B}^{(2)}_3(i,j,k,l)&=&-(\xi^s_{jkl}+\xi^s_{klu}u_j+\xi^s_{luu}u_ju_k+\xi^s_{uuu}u_ju_ku_l)u_{is},\\
\mathcal{B}^{(2)}_4(i,j,k,l)&=&-(\xi^s_{lu}+\xi^s_{uu}u_l)u_{jk}u_{is}.
\end{eqnarray*}
For the third order derivative term $\mathcal{F}_3$, we obtain that
\begin{eqnarray*}
\mathcal{B}^{(3)}_1(i,j,k,l)&=&(\phi_{iu}+\phi_{uu}u_i)u_{jkl}-(\xi^s_{iu}+\xi^s_{uu}u_i)u_su_{jkl},\\
\mathcal{B}^{(3)}_2(i,j,k,l)&=&-(\xi^s_{kl}+\xi^s_{lu}u_k+\xi^s_{uu}u_ku_l+\xi^s_uu_{kl})u_{ijs}-\xi^s_uu_{is}u_{jkl}.
\end{eqnarray*}
For the fourth order derivative term $\mathcal{F}_4$, we derive that
\begin{eqnarray*}
\mathcal{B}^{(4)}_1(i,j,k,l)&=&(\phi_u-\xi^s_uu_s)u_{ijkl},\\
\mathcal{B}^{(4)}_2(i,j,k,l)&=&-(\xi^s_i+\xi^s_uu_i)u_{jkls}.
\end{eqnarray*}
Hence, we have established the fourth order prolongation formula based on the circle summation of the indicators.

\vspace{20pt}

\section{Symmetry group of Monge-Amp\`{e}re equation}\label{sec4}

\noindent In this section, we will classify the symmetry group of Monge-Amp\`{e}re equation \eqref{e1.4}. The steps of Lie's theory are first to find the conditions on infinitesimal generator $\overrightarrow{v}$ of a one parameter symmetry group of the equation. This is a linear condition imposed on $\overrightarrow{v}$ or $pr^{(2)}\overrightarrow{v}$. After comparing the coefficients of like terms on both sides of the equation and solving the resulting complex system of equations, we obtain all infinitesimal generators $\overrightarrow{v}$. Then, using the one-to-one correspondence of Lie algebras with Lie groups, one can recover all symmetry group of the equation. Now, let us give the proof of Theorem \ref{t1.1}.

\noindent\textbf{Proof of Theorem \ref{t1.1}:} Setting the infinitesimal generator of the one parameter symmetry group by
   $$
    \overrightarrow{v}=\xi^i(x,u)\frac{\partial}{\partial x^i}+\phi(x,u)\frac{\partial}{\partial u},
   $$
we have the prolongation formula for $\overrightarrow{v}$ up to second order is given by
  \begin{equation*}
    pr^{(2)}\overrightarrow{v}=\xi^i\frac{\partial}{\partial x^i}+\phi\frac{\partial}{\partial u}+\phi^i\frac{\partial}{\partial u_i}+\phi^{ij}\frac{\partial}{\partial u_{ij}},
  \end{equation*}
where the coefficients $\phi^{i}, \phi^{ij}$ are presented in Section \ref{sec3}. Acting the prolongation $pr^{(2)}\overrightarrow{v}$ to the equation \eqref{e1.4}, it yields that
   \begin{equation}\label{e4.1}
   u^{ij}\phi^{ij}=0,
   \end{equation}
where $[u^{ij}]$ stands for the inverse matrix of $[u_{ij}]$. Since $u^{is}u_{sj}=\delta_{ij}\det D^2u=\delta_{ij}$ for any $i,j$, equation \eqref{e4.1} can also be written as
   \begin{equation}\label{e4.2}
     \eta_{ij}u^{ij}+N\phi_u-2\textstyle\sum\limits_s\xi^s_s-(N+2)\xi^s_uu_s=0,
   \end{equation}
where
  $$
     \eta_{ij}=\phi_{ij}+\phi_{iu}u_j+(\phi_{uj}+\phi_{uu}u_j)u_i-\big[\xi^s_{ij}+\xi^s_{iu}u_j+(\xi^s_{uj}+\xi^s_{uu}u_j)u_i\big]u_s.
  $$
Comparing the coefficients on both sides of \eqref{e4.2}, we obtain that
   \begin{eqnarray}\label{e4.3}
      && {\begin{cases}
          \xi^i_u=0,\ \forall i,\\
          N\phi_u=2\textstyle\sum_s\xi^s_s,\\
          \phi_{ij}=\phi_{uu}=0, \ \forall i,j,\\
          \phi_{iu}u_j+\phi_{uj}u_i-\xi^s_{ij}u_s=0,\ \forall i,j.
         \end{cases}}
   \end{eqnarray}
It follows from the third identity of \eqref{e4.3} that $\phi$ take the form
    \begin{equation}\label{e4.4}
    \phi(x,u)=(C_iu+D_i)x^i+cu+d,\ \ C_i,D_i,c,d\in \mathbb{R}.
    \end{equation}
Substituting \eqref{e4.4} into the last identity of \eqref{e4.3}, we derive that
$$C_iu_j+C_ju_i-\xi^s_{ij}u_s=0,~\forall i,j,$$
and it yields that $\xi^i$ can be represent as
   \begin{equation}
   \xi^i(x,u)=C_jx^jx^i+A^i_jx^j+B^i,\ \ A^i_j,B^i\in \mathbb{R}.
    \end{equation}
Then the second identity of \eqref{e4.3} is equivalent to $$\bigg(\frac{N}{2}+1\bigg)C_sx^s=\frac{Nc}{2}-\textstyle\sum\limits_sA^s_s,$$
and we immediately obtain that $C_i=0$ and $\textstyle\sum_sA^s_s=\frac{Nc}{2}$. Hence, one concludes that
    \begin{equation}
      \begin{cases}
        \xi^i(x,u)=A^i_jx^j+B^i,\\
        \phi(x,u)=D_ix^i+cu+d.
      \end{cases}
    \end{equation}
And the Lie algebra of infinitesimal generators are spanned by the vector fields
  \begin{eqnarray*}
    \overrightarrow{v}^i_1&=&\partial_{x^i},\\
    \overrightarrow{v}_2&=&\partial_u,\\
    \overrightarrow{v}^i_3&=&x^i\partial_u,\\
    \overrightarrow{v}_4&=&a^i_jx^j\partial_{x^i},\\
    \overrightarrow{v}^i_5&=&Nx^i\partial_{x^i}+2u\partial_u,
  \end{eqnarray*}
where $i=1,2,\cdots,N$ and the constants $a^i_j$ satisfy $\textstyle\sum_i a^i_i=0$. Thus, we have completed the proof.  \hfill $\Box$

\vspace{5pt}

As summarized above, we obtain the following group characterization for \eqref{e1.4}.
\begin{theo}
  The Monge-Amp\`{e}re equation \eqref{e1.4} is invariant under group action
    \begin{equation}
      \left[
        \begin{array}{c}
           x\\
           u
        \end{array}
      \right]\to
      \left[
        \begin{array}{c}
           \widetilde{x}\\
           \widetilde{u}
        \end{array}
      \right]=\left[
        \begin{array}{cc}
            A & 0\\
            D & c
        \end{array}
      \right]\left[
        \begin{array}{c}
            x\\
            u
        \end{array}
      \right]+\left[
        \begin{array}{c}
            B\\
            d
        \end{array}
      \right],
    \end{equation}
 where
     \begin{eqnarray*}
      &A=(A^\alpha_\beta)_{N\times N}, \ B=(B^\alpha)_{N\times 1}, \ D=(D_\beta)_{1\times N},&\\
      &c=(c)_{1\times 1}, \ d=(d)_{1\times 1}&
     \end{eqnarray*}
 are matrices satisfying
    \begin{equation}
      A=(\lambda I)\cdot\overline{A}, \ \ \det(\overline{A})=1, \ \ c=\lambda^2
    \end{equation}
 for positive constant $\lambda$.
\end{theo}

\bigskip
\section{Symmetry group of affine maximal type equation}\label{sec5}

\noindent In this section, we classify the symmetry group of the affine maximal type equation \eqref{e1.2}. Since \eqref{e1.2} is a much more complicated fourth order fully nonlinear equation, we need to introduce some notations for simplicity. For any convex solution $u$ of \eqref{e1.2} with $N\geq1$ and  $\theta>0$, we use lower indices like $u_i\equiv D_iu$ to denote the usual derivatives of $u$, and use upper indices $T^i_j\equiv u^{ia}T_{aj}$ to denote the conjugate ones. For example,
   $$
    u^a_{bc}\equiv u^{ad}u_{bcd}, \ \ u^{ab}_c\equiv u^{ai}u^{bj}u_{ijc},\ \ B^{i,j}\equiv u^{ia}u^{jb}B_{a,b}
   $$
and so on will be used below. Now we prove our main classification theorem.

\noindent\textbf{Proof of Theorem \ref{t1.2}:} Elementary computations show that
   $$
     D_ku^{ij}=-u^{ip}u^{jq}u_{pqk}=-u^{ij}_k,
   $$
and hence the first order derivative of $w$ is given by
    \begin{equation}\label{e5.1}
     w_k=-\theta [\det D^2u]^{-\theta-1}U^{pq}u_{pqk}=-\theta wu^a_{ak}.
    \end{equation}
Differentiating once more on both sides of \eqref{e5.1}, we get
    \begin{equation}\label{e5.2}
      w_{kl}=w\big[\theta^2u^a_{ak}u^b_{bl}-\theta(u^a_{akl}-B_{k,l})\big]
    \end{equation}
for $B_{k,l}\equiv u^{ab}_ku_{abl}$. Taking trace on both sides of \eqref{e5.2} and using \eqref{e1.2}, we obtain that
    \begin{eqnarray}\nonumber\label{e5.3}
       u^{ik}_{ik}&=&u^{ij}u^{kl}u_{ijkl}\\ \nonumber
          &=&\theta u^{ij}u^{k}_{ki}u^l_{lj}+u_{ijk}u^{ijk}\\
          &\equiv&\theta v+z,
    \end{eqnarray}
where $v\equiv u^{a}_{ai}u^{bi}_b$ and $z\equiv u_{ijk}u^{ijk}$ are quasi-norm and canonic norm of third order derivatives of $u$ as usual. By \eqref{e5.2}, the equation \eqref{e1.2} is equivalent to
   \begin{equation}\label{e5.4}
     \theta u^{ab}u^{ij}u^{ck}u_{abc}u_{ijk}-u^{ij}u^{kl}u_{ijkl}+u^{ai}u^{bj}u^{ck}u_{abc}u_{ijk}=0.
   \end{equation}
Denoting the infinitesimal generator of the one parameter symmetry group by
   $$
    \overrightarrow{v}=\xi^i(x,u)\frac{\partial}{\partial x^i}+\phi(x,u)\frac{\partial}{\partial u},
   $$
we have the prolongation formula for $\overrightarrow{v}$ up to fourth order is given by
   $$
    pr^{(4)}\overrightarrow{v}=\xi^i\frac{\partial}{\partial x^i}+\phi\frac{\partial}{\partial u}+\phi^i\frac{\partial}{\partial u_i}+\phi^{ij}\frac{\partial}{\partial u_{ij}}+\phi^{ijk}\frac{\partial}{\partial u_{ijk}}+\phi^{ijkl}\frac{\partial}{\partial u_{ijkl}},
   $$
where the coefficients $\phi^{i}, \phi^{ij}, \phi^{ijk}, \phi^{ijkl}$ are presented in Section \ref{sec3}. Since
   $$
     0=\frac{\partial}{\partial u^{ij}}(u^{ab}u_{bc})=\frac{\partial u^{ab}}{\partial u_{ij}}u_{bc}+u^{ab}\frac{\partial u_{bc}}{\partial u_{ij}},
   $$
we derive that
   \begin{eqnarray}
     &&\frac{\partial u^{ad}}{\partial u_{ij}}=-u^{ab}\frac{\partial u_{bc}}{\partial u_{ij}}u^{cd}=-u^{ai}u^{jd},
   \end{eqnarray}
where $u_{ij}$ are regarded as independent variables. Direct calculation shows that $\overrightarrow{v}$ is an infinitesimal generator of a symmetry group of \eqref{e5.4} if and only if
  $$
    \theta u^{ab}u^{ij}u^{ck}\chi_1-\theta u_{abc}u_{ijk}\chi_2-u^{ij}u^{kl}\phi^{ijkl}+u_{ijkl}\chi_3+u^{ai}u^{bj}u^{ck}\chi_1-u_{abc}u_{ijk}\chi_4=0,
  $$
where
   \begin{eqnarray*}
    \chi_1&=&\phi^{abc}u_{ijk}+u_{abc}\phi^{ijk},\\
    \chi_2&=&(u^{ap}u^{bq}u^{ij}u^{ck}+u^{ab}u^{ip}u^{jq}u^{ck}+u^{ab}u^{ij}u^{cp}u^{kq})\phi^{pq},\\
    \chi_3&=&(u^{ip}u^{jq}u^{kl}+u^{ij}u^{kp}u^{lq})\phi^{pq},\\
    \chi_4&=&(u^{ap}u^{iq}u^{bj}u^{ck}+u^{ai}u^{bp}u^{jq}u^{ck}+u^{ai}u^{bj}u^{cp}u^{kq})\phi^{pq},
   \end{eqnarray*}
which can be simplified by
  \begin{eqnarray}\label{e5.6}
  2(\theta u^{ak}_au^{ij}+u^{ijk})\phi^{ijk}-(2\theta u^a_{ak}u^{ijk}+\theta u^{ai}_au^{bj}_b+3B^{i,j}-2u^{aij}_a)\phi^{ij}-u^{ij}u^{kl}\phi^{ijkl}=0.
  \end{eqnarray}
Let us begin by focusing on the fourth order and second order derivative terms in identity \eqref{e5.6}. For one thing, after substituting \eqref{e5.2} into \eqref{e5.6}, we get that the terms for the fourth order derivative, which does not involve the third order derivative, are
  \begin{eqnarray*}
   &&2u^{aij}_a\phi^{ij}+4(\xi^s_i+\xi^s_uu_i)u^{ij}u^{kl}u_{jkls}=0,
  \end{eqnarray*}
where
    \begin{equation}\label{e5.7}
     u^{ab}_{ab}=u^{abc}u_{abc}+\theta u^{kl}u^a_{ak}u^b_{bl}
    \end{equation}
was changed to lower terms by \eqref{e5.3}. So, we conclude that for any $i,j$, there holds
   $$
      \phi_{ij}+\phi_{iu}u_j+(\phi_{uj}+\phi_{uu}u_j)u_i-\big[\xi^s_{ij}+\xi^s_{iu}u_j+(\xi^s_{uj}+\xi^s_{uu}u_j)u_i\big]u_s=0,
   $$
or equivalent to
\begin{equation}\label{e5.8}
\left\{
\begin{aligned}
          &\phi_{ij}=0 \ (\forall i,j),\ \xi^i_{uu}=0\ (\forall i),\\
          &\phi_{uu}=2\xi^i_{iu}\ (\forall i), \ \xi^i_{uj}=0\ (\forall i\neq j),\\
          &2\phi_{iu}=\xi^i_{ii}\ (\forall i),\ \phi_{uj}=\xi^i_{ij}, \ \xi^i_{jj}=0\ (\forall i\neq j).
\end{aligned}
\right.
\end{equation}
For another, the second order derivative term is given by
$$u^{ij}u^{kl}\Big\{\phi_{uu}\textstyle\sum\limits_{\mathrm{circ}}u_{ij}u_{kl} -\textstyle\sum\limits_{\mathrm{circ}}\xi^s_{lu}u_{jk}u_{is}\Big\}=N(N+2)\phi_{uu}-4(N+2)\textstyle\sum\limits_s\xi^s_{su}=0.$$
Note that the identity in the second line of \eqref{e5.8} implies that $N\phi_{uu}=2\textstyle\sum_s\xi^s_{su}$. Hence, we get that $\phi_{uu}=0$, and then $\xi^i_{iu}=0$ for any $i$. Therefore, we conclude that
   \begin{eqnarray}
      &&{\begin{cases}
          \phi(x,u)=(C_iu+D_i)x^i+cu+d, & C_i,D_i,c,d\in \mathbb{R},\\
          \xi^i(x,u)=P^i u+Q^i(x), & P^i\in{\mathbb{R}},\\
          \xi^i_{ij}=C_j\ (\forall i\neq j),\ \xi^i_{ii}=2C_i\ (\forall i).
       \end{cases}}
   \end{eqnarray}

By applying the conclusions obtained above, we can significantly simplify \eqref{e5.6}. Let us now focus on the third order derivative terms. Recall that the fourth order derivative term $u^{ab}_{ab}$ can be expressed in terms of third order derivatives by \eqref{e5.7}. Consequently, the terms in \eqref{e5.6} that may involve the square of third order derivatives are given by
   $$
     2(\theta u^{ak}_au^{ij}+u^{ijk})\chi_5-(2\theta u^a_{ak}u^{ijk}+\theta u^{ai}_au^{bj}_b+3B^{i,j}-2u^{aij}_a)\chi_6-u^{ij}u^{kl}\chi_7,
   $$
where
   \begin{eqnarray*}
    \chi_5&=&\textstyle\sum\limits_{\mathrm{circ}}\Big\{(\phi_u-\xi^s_uu_s)u_{ijk}-(\xi^s_i+\xi^s_uu_i)u_{jks}\Big\},\\
    \chi_6&=&\textstyle\sum\limits_{\mathrm{circ}}\Big\{\phi_{ju}u_i+\phi_uu_{ij}-(\xi^s_{ij}+\xi^s_uu_{ij})u_s -(\xi^s_i+\xi^s_uu_i)u_{js}\Big\},\\
    \chi_7&=&\textstyle\sum\limits_{\mathrm{circ}}\Big\{(\phi_u-\xi^s_uu_s)u_{ijkl}-(\xi^s_i+\xi^s_uu_i)u_{jkls}\Big\}.
   \end{eqnarray*}
If we denote the terms involving $v\equiv u^{ai}_au^b_{bi}$ and $z\equiv u^{abc}u_{abc}$ by $\mathcal{W}_1$, then
   $$
    \mathcal{W}_1=(\phi_u-\xi^s_uu_s)(2z-3z-\theta v-\theta v-z+2\theta v+2z)+4\textstyle\sum\limits_s\xi^s_s(\theta v+z-\theta v-z)\equiv0.
   $$
Next, the terms containing the square of the third order derivative and do not involve $v$ and $z$ in \eqref{e5.6} are given by
   $$
    2(\theta u^{ak}_au^{ij}+u^{ijk})\chi_8-(2\theta u^a_{ak}u^{ijk}+\theta u^{ai}_au^{bj}_b+3B^{i,j})\chi_9,
   $$
where
   \begin{eqnarray*}
    \chi_8&=&-\textstyle\sum\limits_{\mathrm{circ}}(\xi^s_i+\xi^s_uu_i)u_{jks},\\
    \chi_9&=&\textstyle\sum\limits_{\mathrm{circ}}\Big\{\phi_{iu}u_j-\xi^s_{ij}u_s-(\xi^s_i+\xi^s_uu_i)u_{js}\Big\}.
   \end{eqnarray*}
Letting $\mathcal{W}_2$ be the sum of all terms involving $u^{ak}_au^{b}_{bs}$, we have
   $$
     \mathcal{W}_2=-2\theta u^{ak}_au^b_{bs}\xi^s_uu_k-2\theta u^{ak}_au^b_{bs}\xi^s_k-\theta u^{ai}_au^{bj}_b(2\phi_{iu}u_j-\xi^s_{ij}u_s-2\xi^s_uu_iu_{js} -2\xi^s_iu_{js})=0,
   $$
which is equivalent to $\phi_{iu}=0\ (\forall i)$ and $\xi^s_{ij}=0\ (\forall s,i,j)$. Furthermore, by examining the terms containing $u^{ak}_au^i_{ks}$, denoted by $\mathcal{W}_3$, we obtain that
   $$
     \mathcal{W}_3=-4\theta u^{ak}_au^i_{ks}\xi^s_uu_i-4\theta u^{ak}_au^i_{ks}\xi^s_i-2\theta u^{ak}_au^{ij}_k(-\xi^s_{ij}u_s-2\xi^s_iu_{js}-2\xi^s_uu_iu_{js})=0,
   $$
which implies ${\xi^s_{ij}=0}\ (\forall s,i,j)$. By checking the terms containing $B^{i,j}$, denoted by $\mathcal{W}_4$, we derive that
   $$
     \mathcal{W}_4=-6B^{i,}_{,s}\xi^s_uu_i-6B^{i,}_{,s}\xi^s_i-3B^{i,j}(-2\xi^s_iu_{js}-2\xi^s_uu_iu_{js})\equiv0.
   $$
Finally, denoting the terms containing $u_{ijk}$ by $\mathcal{W}_5$, we get that
   $$
     \mathcal{W}_5=-2\theta u^{ak}_au^{ij}\xi^s_u\chi_{10}-6u^{ijk}\xi^s_uu_{ij}u_{ks}+u^{ij}u^{kl}\xi^s_u\chi_{11}=0,
   $$
where
   \begin{eqnarray*}
    \chi_{10}&=&u_{ij}u_{ks}+u_{jk}u_{is}+u_{ik}u_{js},\\
    \chi_{11}&=&u_{ijk}u_{ls}+3u_{ijl}u_{ks}+2u_{ij}u_{kls}+4u_{jk}u_{ils},
   \end{eqnarray*}
which is equivalent to
    \begin{equation*}
      -2(N+2)\bigg(\theta-\frac{N+1}{N+2}\bigg)u^b_{bs}\xi^s_u=0\Leftrightarrow\begin{cases}
         \xi^i_u=0, & \mbox{if } \theta\not=\frac{N+1}{N+2},\\
         0=0, & \mbox{if } \theta=\frac{N+1}{N+2}.
      \end{cases}
    \end{equation*}

Therefore, we conclude that for $\theta=\frac{N+1}{N+2}$,
   \begin{equation}
     \begin{cases}
      \xi^i(x,u)=P^i u+Q^i_j x^j+R^i, & P^i, Q^i_j, R^i\in{\mathbb{R}},\\
      \phi(x,u)=D_ix^i+cu+d, & D_i,c, d\in{\mathbb{R}}
     \end{cases}
   \end{equation}
and for $\theta\not=\frac{N+1}{N+2}$,
   \begin{equation}
      \begin{cases}
       \xi^i(x,u)=Q^i_j x^j+R^i, &  Q^i_j, R^i\in{\mathbb{R}},\\
       \phi(x,u)=D_ix^i+cu+d, & D_i,c, d\in{\mathbb{R}}.
     \end{cases}
   \end{equation}
Hence, the Lie algebra of infinitesimal generators are spanned by the vector fields
    \begin{eqnarray*}
    \overrightarrow{v}_1^i&=&\partial_{x^i},\\
    \overrightarrow{v}_2&=&\partial_u,\\
    \overrightarrow{v}_3&=&u\partial_u,\\
    \overrightarrow{v}_4^i&=&x^i\partial_u,\\
    \overrightarrow{v}_5^{ij}&=&x^j\partial_{x^i},\\
    \overrightarrow{v}_6^i&=&u\partial_{x^i}
  \end{eqnarray*}
for $\theta=\frac{N+1}{N+2}$, and $\overrightarrow{v}_1^i$,$\overrightarrow{v}_2$,$\overrightarrow{v}_3$,$\overrightarrow{v}_4^i$,$\overrightarrow{v}_5^{ij}$ for $\theta\not=\frac{N+1}{N+2}$. Thus, the proof was done.  \hfill $\Box$

\vspace{5pt}

\begin{rema}
Since each solution to the Monge-Amp\`{e}re equation \eqref{e1.4} also satisfies the affine maximal type equation \eqref{e1.2}, it follows that the symmetry group of \eqref{e1.2} is larger than that of \eqref{e1.4}.
\end{rema}

Summing as above, we obtain the following group characterization for \eqref{e1.2}.
\begin{theo}
  For $\theta=\frac{N+1}{N+2}$, the affine maximal type equation \eqref{e1.2} is invariant under group action
    \begin{equation}\label{e5.12}
      \left[
        \begin{array}{c}
           x\\
           u
        \end{array}
      \right]\to
      \left[
        \begin{array}{c}
           \widetilde{x}\\
           \widetilde{u}
        \end{array}
      \right]=\left[
        \begin{array}{cc}
            Q & P \\
            D & c
        \end{array}
      \right]\left[
        \begin{array}{c}
            x\\
            u
        \end{array}
      \right]+\left[
        \begin{array}{c}
            R\\
            d
        \end{array}
      \right],
    \end{equation}
  where
     \begin{eqnarray*}
      &Q=(Q^\alpha_\beta)_{N\times N}, \ P=(P^\alpha)_{N\times 1}, \ R=(R^\alpha)_{N\times1},&\\
      &D=(D_\beta)_{1\times N}, \ c=(c)_{1\times1}, \ d=(d)_{1\times1}&
     \end{eqnarray*}
  are matrices such that
     $$
      \left[
       \begin{array}{cc}
         Q & P\\
         D & c
       \end{array}
       \right]\in GL(N+1).
     $$
   More precisely, the group action is global for $P=0$ and only local for $P\neq 0$. For $\theta\not=\frac{N+1}{N+2}$, there are only group actions \eqref{e5.12} for $P=0$.
\end{theo}

\vspace{20pt}

{\small\noindent\textbf{Acknowledgments}\enspace The author (SZ) would like to express his deepest gratitude to Professors Xi-Ping Zhu, Kai-Seng Chou, Xu-Jia Wang and Neil Trudinger for their constant encouragements and warm-hearted helps. This paper is also dedicated to the memory of Professor Dong-Gao Deng.}


\vspace{20pt}

\begin{thebibliography}{}

\bibitem{A} M. Abreu, {\it K\"{a}hler geometry of toric varieties and extremal metrics}, Int. J. Math., {\bf9} (1998), 641-651.\\

\bibitem{Ca} E. Calabi, {\it Hypersurfaces with maximal affinely invariant area}, Amer. J. Math., {\bf104} (1982), 91-126.\\

\bibitem{Ch} S.S. Chern, {\it Affine minimal hypersurfaces, Minimal submanifolds and geodesics}, Proc. Japan-United States Sem., Tokyo, 1977, 17-30.\\

\bibitem{CL} K.S. Chou and G.X. Li, {\it Optimal systems and invariant solutions for the curve shortening problem}, Comm Anal. Geom., {\bf10} (2002), 241-274.\\

\bibitem{CM} P.A. Clarkson and E.L. Mansfield, \emph{Symmetry reductions and exact solutions of a class of nonlinear heat equations}, Phys. D, {\bf70} (1994), 250-288.\\

\bibitem{CQ} K.S. Chou and C.Z. Qu, {\it Optimal systems and group classification of $(1+2)$-dimensional heat equation}, Acta Appl. Math., {\bf83} (2004), 257-287.\\

\bibitem{CW1} K.S. Chou and X.J. Wang, {\it A variational theory of the Hessian equation}, Comm. Pure Appl. Math., {\bf54} (2002), 1029-1064.\\

\bibitem{CDG} X.P. Chen, S.Z. Du and T.P. Guo, {\it The Liouville theorem of a torsion system and its application to the symmetry group of a porous medium type equation on symmetric spaces}, J. Lie theory, {\bf31} (2021), 393-411.\\

\bibitem{CNS} L. Caffarelli, L. Nirenberg and J. Spruck, {\it The Dirichlet problem for nonlinear second-order elliptic equations. I. Monge-Amp\`{e}re equation}, Comm. Pure Appl. Math., {\bf37} (1984), 369-402.\\

\bibitem{Do} S.K. Donaldson, {\it Interior estimates for solutions of Abreu's equation}, Collect. Math., {\bf56} (2005), 103-142.\\

\bibitem{Du1} S.Z. Du, {\it Bernstein problem of affine maximal type hypersurfaces on dimension $N\geq3$}, J. Differential Equations, {\bf269} (2020), 7429-7469.\\

\bibitem{Du2} S.Z. Du, {\it Non-quadratic Euclidean complete affine maximal type hypersurfaces for $\theta\in(0,(N-1)/N]$}, J. Geom. Anal., {\bf229} (2024), 1-19.\\

\bibitem{JL} F. Jia and A.M. Li, {\it A Bernstein property of affine maximal hypersurfaces}, Ann. Global. Anal. Geom., {\bf23} (2003), 359-372.\\

\bibitem{JL2} F. Jia and A.M. Li, {\it A Bernstein property of some fourth order partial differential equations}, Results Math., {\bf56} (2009), 109-139.\\

\bibitem{L3} P.L. Lions, N.S. Trudinger and J.I.E. Urbas, {\it The Neumann problem for equations of Monge-Amp\`{e}re type}, Comm. Pure Appl. Math., {\bf39} (1986), 539-563.\\

\bibitem{M} J.A. McCoy, {\it A Bernstein property of solutions to a class of prescribed affine mean curvature equations}, Ann. Global. Anal. Geom., {\bf32} (2007), 147-165.\\

\bibitem{O} P.J. Olver, {\it Applications of Lie groups to differential equations}, second edition, Graduate Texts in Mathematics, 107. Springer-Verlag, New York, 1993.\\

\bibitem{OST} P.J. Olver, G. Saporo and A. Tannenbaum, {\it Classification and uniqueness of invariant geometric flows}, C. R. Acad. Sci. Paris Ser. I Math., {\bf319} (1994), 339-344.\\

\bibitem{P} A.V. Pogorelov, {\it On the improper affine hyperspheres}, Geom. Dedicata, {\bf1} (1972), 33-46.\\

\bibitem{T} N.S. Trudinger, {\it On the Dirichlet problem for Hessian equations}, Acta Math., {\bf175} (1995), 151-164.\\

\bibitem{Tso} K. Tso, {\it On a real Monge-Amp\`{e}re functional}, Invent. Math., {\bf101} (1990), 425-448.\\

\bibitem{TW1} N.S. Trudinger and X.J. Wang, {\it The Bernstein problem for affine maximal hypersurfaces}, Invent. Math., {\bf140} (2000), 399-422.\\

\bibitem{TW2} N.S. Trudinger and X.J. Wang, {\it Affine complete locally convex hypersurfaces}, Invent. Math., {\bf150} (2002), 45-60.\\

\bibitem{TW3} N.S. Trudinger and X.J. Wang, {\it Boundary regularity for the Monge-Amp\`{e}re and affine maximal surface equations}, Ann. Math., {\bf167} (2008), 993-1028.\\

\bibitem{W} X.J. Wang, {\it A class of fully nonlinear elliptic equations and related functionals}, Indiana Univ. Math. J., {\bf43} (1994), 25-54.\\

\bibitem{WYW} W.F. Wo, S.X. Yang and X.L. Wang, {\it Group invariant solutions to a centro-affine invariant flow}, Arch. Math. (Basel), {\bf108} (2017), 33-46.\\

\bibitem{Z} B. Zhou, {\it The Bernstein theorem for a class of fourth order equations}, Calc. Var. Partial Differential Equations, {\bf43} (2012), 25-44.



\end{thebibliography}
\end{document}